\documentclass{amsart}
\usepackage{amsmath,amsthm,amsfonts,amscd,eucal,latexsym}
\usepackage{graphicx}
\usepackage{epsfig}
\usepackage{subfigure}

\numberwithin{equation}{section}

\newcommand{\cb}{{\mathcal B}}
\newcommand{\ct}{{\mathcal T}}
\newcommand{\cu}{{\mathcal U}}

\newcommand{\bn}{{\mathbb N}}
\newcommand{\br}{{\mathbb R}}

\renewcommand{\a}{\alpha}
\renewcommand{\d}{\delta}      
\newcommand{\eps}{\varepsilon} 
\renewcommand{\k}{\kappa}
\renewcommand{\l}{\lambda}   
\newcommand{\m}{\mu}
\newcommand{\n}{\nu}
\renewcommand{\r}{\rho}
\newcommand{\f}{\varphi}
\renewcommand{\c}{\chi}

\newcommand{\subc}{\underline{\delta}}
\newcommand{\supc}{\overline{\delta}}
\newcommand{\subd}{\underline{d}}
\newcommand{\supd}{\overline{d}}

\newcommand{\ov}{\overline}
\newcommand{\itm}[1]{\item{$(#1)$}}

\newcommand{\e}[1]{\text{e}^{#1}}
\newcommand{\subg}{\underline{g}}
\newcommand{\supg}{\overline{g}}
\newcommand{\B}{\overline{B}}
\newcommand{\floor}[1]{\lfloor #1 \rfloor}

\newcommand{\convGH}{\begin{smallmatrix} \\ \longrightarrow \\
\text{GH}\end{smallmatrix}}

\newcommand{\convPGH}{\begin{smallmatrix} \\ \longrightarrow \\
\text{pGH}\end{smallmatrix}}

\newtheorem{thm}{Theorem}[section]

\newtheorem{prop}[thm]{Proposition}
\newtheorem{lem}[thm]{Lemma}
\theoremstyle{definition}
\newtheorem{dfn}[thm]{Definition}

\theoremstyle{remark}
\newtheorem{rmk}[thm]{Remark} 
 
\newtheorem{Assump}[thm]{Assumption}
\begin{document}
\title{Tangential dimensions I. Metric spaces}
\author{Daniele Guido, Tommaso Isola}
\address{
Dipartimento di Matematica, Universit\`a di Roma ``Tor
Vergata'', I--00133 Roma, Italy.}
\email{guido@mat.uniroma2.it, isola@mat.uniroma2.it}
\date{}

\begin{abstract}
    Pointwise tangential dimensions are introduced for metric spaces. 
    Under regularity conditions, the upper, resp.  lower, tangential
    dimensions of $X$ at $x$ can be defined as the supremum, resp. 
    infimum, of box dimensions of the tangent sets, {\it a la Gromov},
    of $X$ at $x$.  Our main purpose is that of introducing a tool
    which is very sensitive to the ``multifractal behaviour at a
    point" of a set, namely which is able to detect the
    ``oscillations" of the dimension at a given point.  In particular
    we exhibit examples where upper and lower tangential dimensions
    differ, even when the local upper and lower box dimensions
    coincide.  Tangential dimensions can be considered as the
    classical analogue of the tangential dimensions for spectral
    triples introduced in \cite{GuIs9}, in the framework of Alain
    Connes' noncommutative geometry \cite{Co}.
\end{abstract}

\maketitle

 \setcounter{section}{0}

 \section{Introduction.}\label{sec:zeroth}
 
 Dimensions can be seen as a tool for measuring the non-regularity, or
 fractality, of a given object.  Non-integrality of the dimension is a
 first sign of non-regularity.  A second kind of non-regularity is
 related to the fact that the dimension is not a global constant. 
 This may happen in two ways: either the dimension varies from point
 to point, or it has an oscillating behavior at a point. 
 Indeed dimensions are often defined as limits, and an oscillating
 behavior means that the upper and lower versions of the considered
 dimension are different.  Our main goal here is to introduce a local
 dimension that is able to maximally detect such an oscillating
 behavior, namely for which the upper and lower determinations form a
 maximal dimensional interval.  With this aim, we shall define the
 upper and lower tangential dimension for a metric space.
 We mention at this point that such dimensions, which are
 presented here in a completely "classic" way, have been introduced
 first for noncommutative spaces \cite{GuIs9}, where their definition
 is purely noncommutative, depending on the oscillating behavior of
 the eigenvalues of the Dirac operator, which may imply that the
 (singular) traceability exponents form an interval, rather than a singleton.

 The name tangential is motivated here by the fact that, under
 suitable hypotheses, such dimensions are the supremum, resp. 
 infimum, of the local dimensions of the tangent sets for the given
 space.  The notion of tangent set (or rather tangent cone, cf. 
 Remark \ref{rem:tangenti}) for a metric space is due to Gromov
 \cite{Gromovbook}.  A tangent set of a metric space $X$ at a point
 $x$ is any limit point of the family of its dilations, for the
 dilation parameter going to infinity, taken in the pointed
 Gromov-Hausdorff topology.

 \begin{figure}[ht]
     \centering
     \subfigure{
     \psfig{file=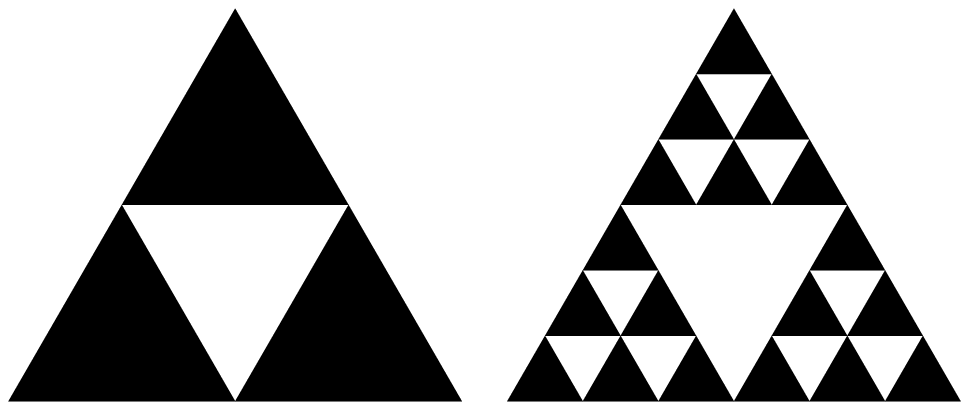,height=1.5in}}
     \hspace{0.3 in}
     \subfigure{
     \psfig{file=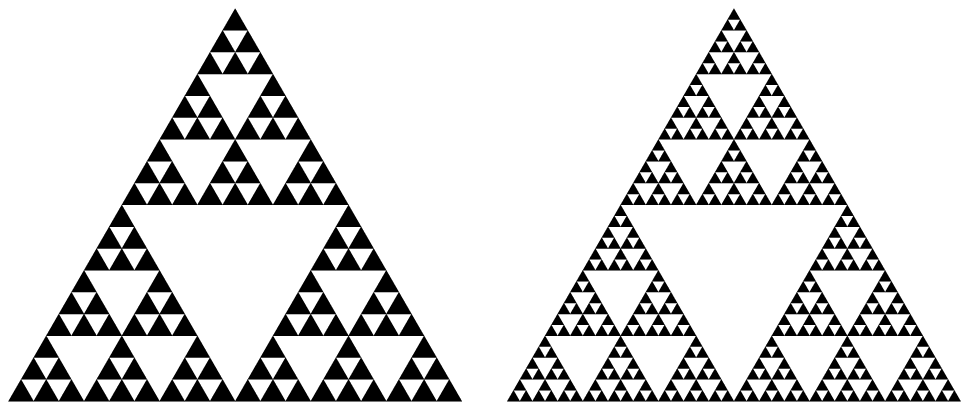,height=1.5in}}
     \caption{Modified Sierpinski}
     \label{fig:Sierp}
 \end{figure}

 As an example we mention some fractals considered in \cite{Hambly}. 
 They are constructed as follows.  At each step the sides of an
 equilateral triangle are divided in $q\in\bn$ equal parts, so as to
 obtain $q^{2}$ equal equilateral triangles, and then all downward
 pointing triangles are removed, so that $\frac{q(q+1)}{2}$ triangles
 are left.  Setting $q_{j}=2$ if $(k-1)(2k-1)<j\leq (2k-1)k$ and
 $q_{j}=3$ if $ k(2k-1)<j\leq k(2k+1)$, $k=1,2,\dots$, we get a
 translation fractal with dimensions given by \cite{GuIs11b}
 $$
 \subc=\frac{\log3}{\log2}<\subd=\supd=\frac{\log18}{\log6}
 <\supc=\frac{\log6}{\log3},
 $$
 where $\subc,\supc,\subd,\supd$ denote the lower tangential, the
 upper tangential, the lower local and the upper local dimensions. 
 The first four steps ($q=2,3,3,2$) of the procedure above are shown
 in Figure \ref{fig:Sierp}.

 \begin{figure}[ht]
     \centering
     \subfigure{
     \psfig{file=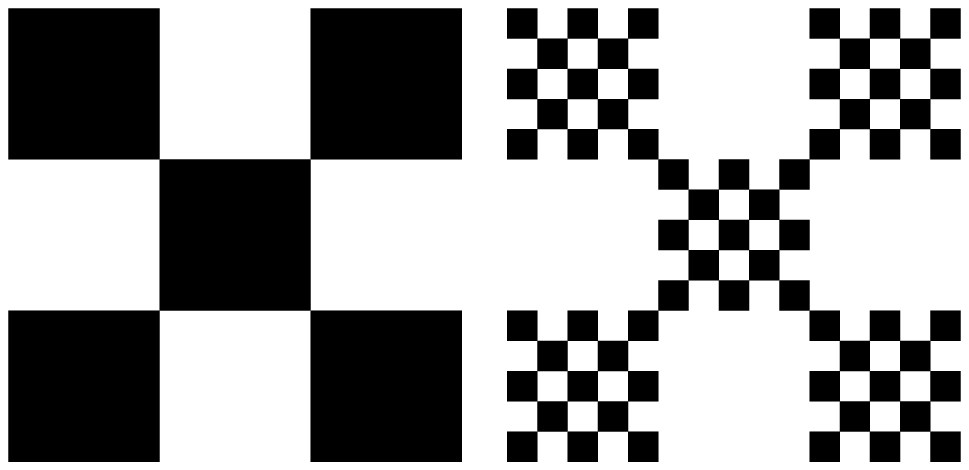,height=1.5in}}
     \hspace{0.3 in}
     \subfigure{
     \psfig{file=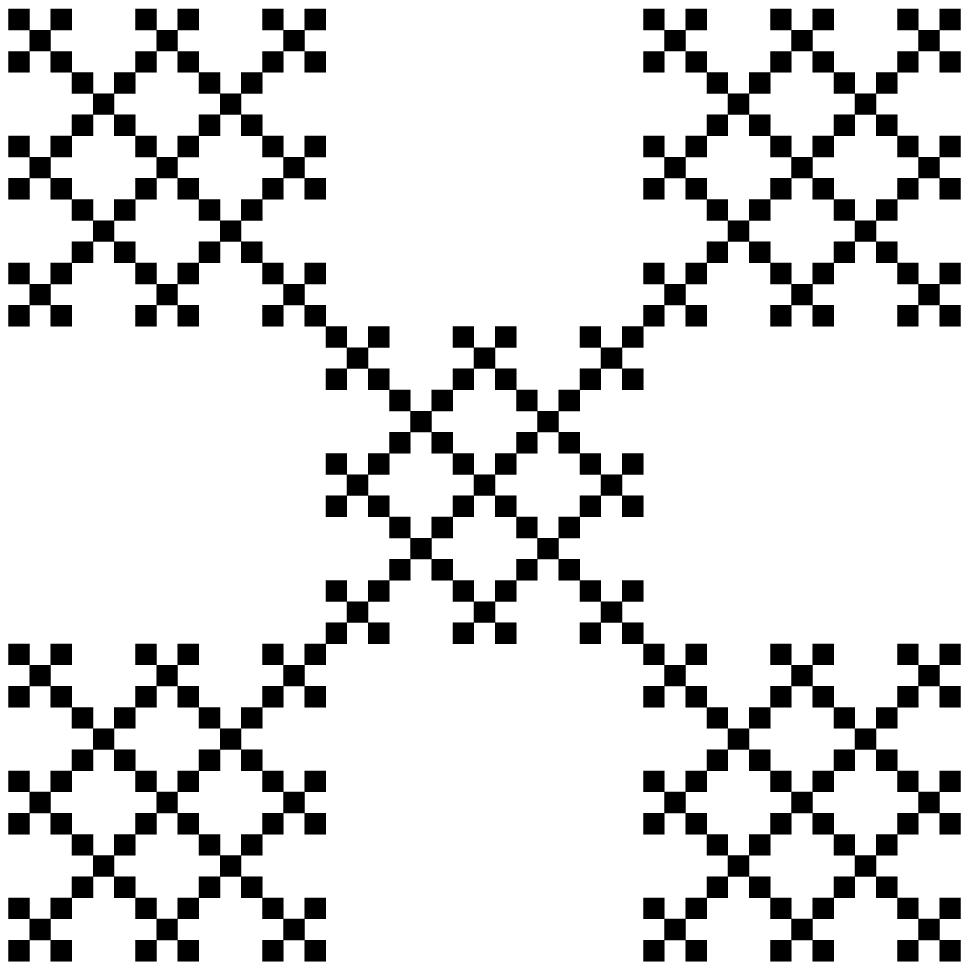,height=1.5in}}
     \caption{Modified Vicsek}
     \label{fig:Vics}
 \end{figure}

 The procedure considered above can, of course, be applied also to
 other shapes.  For example, at each step the sides of a square are
 divided in $2q+1$, $q\in\bn$, equal parts, so as to obtain
 $(2q+1)^{2}$ equal squares, and then $2q(q+1)$ squares are removed,
 so that to remain with a chessboard.  In particular, we may set
 $q_{j}=2$ if $k(2k+1)<j\leq (2k+1)(k+1)$ and $q_{j}=1$ if $
 k(2k-1)<j\leq k(2k+1)$, $k=0,1,2,\dots$, getting a translation
 fractal with dimensions given by \cite{GuIs11b}
 $$
 \subc=\frac{\log5}{\log3}<\subd=\supd=\frac{\log65}{\log15}
 <\supc=\frac{\log13}{\log5}.
 $$

 The first three steps ($q=1,2,1$) of this procedure are shown in Figure
 \ref{fig:Vics}.
 
 The fractals considered above belong to a general class of fractals,
 called translation fractals.  In a forthcoming paper \cite{GuIs11b},
 we shall show that for such fractals the metric tangential dimensions
 coincide with the tangential dimensions of an invariant measure, in
 this way obtaining an explicit formula for the dimensions.
 
 Translation fractals can also be
 studied from a noncommutative point of view, and commutative and
 noncommutative tangential dimensions coincide.  This follows for
 translation fractals in $\br$ simply comparing the formulas given in
 \cite{GuIs9} and those given in \cite{GuIs11b}.  The analysis of translation
 fractals in $\br^{n}$ and their tangential dimensions from a
 noncommutative point of view is contained in \cite{GuIs10}.

\section{Tangent sets of a metric space}

Tangent sets of metric spaces at a point have been defined by Gromov, 
cf. \cite{Gromovbook,BHbook}.

If $(X,d)$ is a metric space, we shall denote by $B(x,r)$ the open
ball $\{y\in X:d(x,y)<r\}$, by $\ov{B}(x,r)$ the closed ball
$\{y\in X:d(x,y)\leq r\}$ and by $\ov{B(x,r)}$ the closure of $B(x,r)$; 
moreover $B_{\eps}(E) := \{x\in X: \inf_{y\in E} d(x,y) <\eps\}$, 
for $E\subset X$.

Let us recall that the Gromov-Hausdorff distance $d_{GH}(X,Y)$ between
two metric spaces $X$ and $Y$ is defined as the infimum of the
$\eps>0$ such that there are isometric embeddings $\f_{X}$, $\f_{Y}$
of $X$ and $Y$ into a metric space $Z$ for which $\f_{X}(X)\subset
B_{\eps}(\f_{Y}(Y))$ and $\f_{Y}(Y)\subset B_{\eps}(\f_{X}(X))$.  This
is indeed a distance between isometry classes of compact metric
spaces.

In case of noncompact (proper) metric spaces one considers the pointed
Gromov-Hausdorff topology, which can be equivalently defined as

\begin{itemize}
    \item[$(1)$] a neighbourhood base consists of the sets $U^{\eps}(X,x)$,
    $(X,x)$ a pointed metric space, $\eps\in(0,\frac12)$, where
    $U^{\eps}(X,x) := \{ (Y,y) : d_{pGH}((X,x),(Y,y))<\eps\}$, and
    $d_{pGH}((X,x),(Y,y))$ is the infimum of the $\eps>0$ for which
    there is a compatible metric $d$ on the disjoint union of $X$ and
    $Y$ s.t. $d(x,y)<\eps$, $\B_{X}(x,\frac{1}{\eps}) \subset
    B_{\eps}(Y)$, $\B_{Y}(y,\frac{1}{\eps}) \subset B_{\eps}(X)$.
    
    \item[$(2)$] a neighbourhood base consists of the sets
    $V^{R,\eps}(X,x)$, $(X,x)$ a pointed metric space,
    $R>0,\,\eps\in(0,1)$, where $V^{R,\eps}(X,x) := \{ (Y,y) :
    d^{R}((X,x),(Y,y))$ $<$ $\eps\}$, and $d^{R}\left((X,x),(Y,y)\right)$ is
    defined as the infimum of the $\eps>0$ such that there are
    isometric embeddings $\f_{X}$, $\f_{Y}$ of $X$ and $Y$ into a
    metric space $(Z,d)$ for which $d(\f_{X}(x),\f_{Y}(y))<\eps$,
    $\f_{X}(\B_{X}(x,R)) \subset B_{\eps}(\f_{Y}(Y))$ and
    $\f_{Y}(\B_{Y}(y,R)) \subset B_{\eps}(\f_{X}(X))$.
\end{itemize}

 On the isometry classes of proper metric spaces it is a Hausdorff
 topology.  Since this topology is separable, it is determined by its
 converging sequences; indeed it is equivalently defined as follows.

 \begin{prop} {\rm \cite{Gromovbook}} \label{Prop:Gromovbook}
      $(X_{n},x_{n})$ converges to $(X,x)$ in the pointed
      Gromov-\-Hausdorff topology if and only if, for any $R>0$ there
      exists a positive infinitesimal sequence $\eps_{n}$ such that,
      for any $\eta>0$ there is $n_{0}\in\bn$ such that for all
      $n>n_{0}$ there are isometric embeddings $\f_{n}$, $\f$ of
      $\ov{B}_{X_{n}}(x_{n},R+\eps_{n})$ and $\ov{B}_{X}(x,R)$ into a
      metric space $(Z_{n},d_{n})$ for which
      $d_{n}(\f_{n}(x_{n}),\f(x))<\eta$,
      $\f_{n}(\ov{B}_{X_{n}}(x_{n},R+\eps_{n})) \subset
      B_{\eta}(\f(\B_{X}(x,R))$ and $\f(\B_{X}(x,R))) \subset
      B_{\eta}(\f_{n}(\B_{X_{n}}(x,R+\eps_{n})))$.
 \end{prop}
 
    From the previous characterization one easily gets

\begin{prop}\label{Prop:conseq} 
    \item[$(i)$] If $(X_{n},x_{n})$ converge to $(X,x)$ in the
    pointed Gromov-Hausdorff topology, then, possibly passing to a
    subsequence, $\ov{B_{X_{n}}(x_{n},R)} \convGH B$, with
    $\ov{B_{X}(x,R)} \subseteq B\subseteq \ov{B}_{X}(x,R)$.

    \item[$(ii)$] If $X_{n}$ is an increasing sequence of proper 
    spaces such that the completion $X$ of $\cup_{n} X_{n}$ is 
    proper, then $(X_{n},x)\convPGH (X,x)$ for any $x\in X_{1}$.
 \end{prop}

\begin{dfn}
    Let $(X,d)$ be a metric space, $x\in X$.  A {\it tangent set} of
    $X$ at $x$ is any limit point, for $t\to\infty$, of $(X,x,td)$ in
    the pointed Gromov-Hausdorff topology, where $td$ denotes the
    rescaled distance by the parameter $t$.  We write also $tX$ for
    $(X,x,td)$ when the metric and $x$ are clear from the context.  We
    shall denote by $\ct_{x} X$, and call the {\it tangent cone} of
    $X$ at $x$, the family of tangent sets of $X$ at $x$.  A {\it
    tangent ball} of $X$ at $x$ is any ball centered in $x$ of some
    tangent set $T\in\ct_{x} X$.
 \end{dfn}

 \begin{prop}
     Let $(X,x)$ be such that 
     \begin{equation}\label{cptcond}
     \limsup_{r\to0}n_{\l r}(\ov{B}_{X}(x,r))<\infty\quad \forall \l>0.
     \end{equation}
     Then $\ct_{x} X$ is not empty.  Indeed, given any sequence
     $t_{n}\to+\infty$, there exists a subsequence $t_{n_{k}}$ for
     which $(X,x,t_{n_{k}}d)$ converges to a unique proper space in
     the pointed Gromov-Hausdorff topology.
 \end{prop}
 
 \begin{proof}
     Follows from the Gromov compactness criterion \cite{Gromovbook}.
 \end{proof}

 \begin{rmk}\label{rem:tangenti}
     \itm{i} A tangent set cannot be empty, since it necessarily
     contains $x$.  It may happen that $\ct_{x} X$ is empty, namely
     that $(X,x,td)$ has no limit points.  
     
     \itm{ii} If $X$ is a manifold, the tangent set at $x$ is unique
     and coincides with the ordinary tangent space (cf. 
     \cite{Gromovbook}).
     
     \itm{iii} $\ct_{x} X$ is indeed a cone, namely it is dilation
     invariant.  In fact, if $(T,d_{T})$ is a tangent set of $X$ at
     $x$ given by the converging sequence $(X,x,t_{n}d)$, and $\a>0$,
     then $(X,x,\a t_{n}d)$, converges to $(T,\a d_{T})$.  As a
     consequence, if $\ct_{x} X$ consists of a unique set, such set is
     a cone.  Since this case has been usually considered, one usually
     refers to such metric space with the name of Gromov tangent cone.
     
     \itm{iv} If all the metric spaces $X_{n}$ are subsets of the same 
     proper metric space $Z$, pointed Gromov-Hausdorff convergence may be 
     replaced by the Attouch-Wets convergence, \cite{Beer}. Let us note 
     that in this case we do not need to specify a point in $Z$.

     \itm{v} If the ambient space $Z$ is dilation invariant, e.g.
     $Z=\br^{n}$, then the dilations of a given subset are still
     subsets of $Z$, hence the tangent sets can be defined as
     Attouch-Wets limits, and are subsets of $Z$, as in \cite{BeFi}. 
     Even if the two topologies do not coincide, the families of
     tangent sets at a given point do.
 \end{rmk}

 We conclude this section by computing explicitly the tangent cone of
 a self-similar fractal at the points which are invariant for some of
 the dilations generating the fractal.

\begin{thm}\label{selfsimilar}
    Let $F$ be a self-similar fractal, $w$ one of the generating 
    similarities,  $x=wx$. 
    The tangent cone $\ct_{x}F$ consists exactly of all dilations of 
    $Z=\overline{\bigcup_{n\in\bn}w^{-n}F}$.
\end{thm}

\begin{proof}
    Let us observe that, given $r_{n}\to +\infty$, $\l>1$, we may find
    $n_{k}, m_{k}\in \bn$, $c>0$ such that $c\l^{m_{k}} / r_{n_{k}}
    \to1$.  Hence, if $T$ is a tangent set of $F$ at $x$, with
    $r_{n}F\convPGH T$, we have $\l^{m_{k}}F\convPGH (1/c) T$.  
    Therefore it is enough to show that $\l^{m}F\convPGH Z$.
    Indeed $(\l^{m}F,x)$ is isometric to $(w^{-n}F,x)$ as pointed 
    metric spaces, and $w^{-n}F$ is an increasing sequence of proper 
    metric spaces, therefore the thesis follows from Proposition 
    \ref{Prop:conseq} (ii).
\end{proof}

Let $\vec{q}=\{q_{j}\}$ be a sequence of natural numbers, and $S(\vec{q})$
be the corresponding fractal constructed as in Fig.  \ref{fig:Sierp}. 
Let us observe that if $q_{j}\equiv q$, then we get the $q$-Sierpinski
triangle $S(q)$.

\begin{thm}
    If for any $p\in\bn$ there is a $\ov{j}$ such that $q_{j}=q$ for
    $\ov{j}\leq j\leq \ov{j}+p$, then the tangent cone of
    $S(\vec{q})$ at one of its extreme points contains the tangent
    space of $S(q)$ at one of its extreme points.
\end{thm}

\begin{proof}
    Clearly for any $p$ we may indeed find an increasing sequence
    $j_{n}$ such that $q_{j}=q$ for $j_{n}-p\leq j\leq j_{n}+p$.  If
    we set $Q_{j}=\prod_{i=1}^{j}q_{i}$, then, for any $r< q^{p}$, the
    Hausdorff distance between the ball of radius $r$ of $Q_{j_{n}}
    S(\vec{q})$ and the ball of radius $r$ of $Q_{j_{n}} S(q)$ is less
    than $q^{-p}$, which implies that the pointed Gromov-Hausdorff
    limit of $Q_{j_{n}} S(\vec{q})$ coincides with the pointed
    Gromov-Hausdorff limit of $Q_{j_{n}} S(q)$.  
\end{proof}

\begin{rmk}
    While in the first example, the tangent sets are
    described by one (dilation) parameter, in the second example a 
    second parameter $q$ appears.  In a sense this
    shows that the higher is the regularity of the set (around
    the point $x$), the smaller is its tangent cone.  In the case of
    the Sierpinski triangle $S$, the explicit description of the
    tangent set to a point $x$ can be extended easily to all points
    which are obtained by applying a product of similarities to one of
    the three extremal points of $S$.
\end{rmk}

\section{Tangential dimensions}

 \subsection{Definition of tangential dimensions and connection with 
 tangent sets}
 
 Let $(X,d)$ be a metric space, $E\subset X$.  Let us denote by
 $n(r,E) \equiv n_{r}(E)$, resp.  $\ov{n}(r,E) \equiv
 \overline{n}_{r}(E)$, the minimum number of open, resp.  closed,
 balls of radius $r$ necessary to cover $E$, and by $\n(r,E)
 \equiv \n_{r}(E)$ the maximum number of disjoint open balls of $E$ of radius
 $r$ contained in $E$.
 
 \begin{dfn}
     Let $(X,d)$ be a metric space, $E\subset X$, $x\in E$.  We call
     {\it upper, resp.  lower tangential dimension} of $E$ at $x$ the
     (possibly infinite) numbers
     \begin{align*}
	\underline{\d}_{E}(x) & := \liminf_{\l \to 0} \liminf_{r \to
	0} \frac{\log n(\l r, E\cap\B(x,r))}{\log 1/\l}, \\
	\overline{\d}_{E}(x) & := \limsup_{\l \to 0} \limsup_{r \to
	0} \frac{\log n(\l r, E\cap\B(x,r))}{\log 1/\l}.
    \end{align*}
 \end{dfn}
 
  \begin{prop}\label{equiv-form}
     Nothing changes in the previous definition if one replaces $n$
     with $\n$ or with $\overline{n}$, or $E\cap\B(x,r)$ with $E \cap
     B(x,r)$.  Moreover, if $E$ is closed in $X$, one can replace
     $E\cap\B(x,r)$ also with $\ov{E \cap B(x,r)}$.
 \end{prop}
 \begin{proof}
     The statements about $\n$ and $\ov{n}$ follow from (see e.g. \cite{GuIs6})
     \begin{align}
	 n_{2r}(E)&\leq \n_{r}(E) \leq n_{r}(E)\label{n-ni}\\
	 n_{2r}(E)&\leq \ov{n}_{r}(E) \leq n_{r}(E).\label{n=ovn}
     \end{align}
     From $B(x,r) \subset \B(x,r) \subset B(x,2r)$, and $E\cap B(x,r)
     \subset \ov{E\cap B(x,r)} \subset E \cap\B(x,r)$, if $E$ is
     closed, follow the other statements.
 \end{proof}

 We want to give a geometric interpretation of the (lower and upper)
 tangential dimensions.  We need some auxiliary results.

\begin{prop}\label{semicontBis}
    \item{$(i)$} For any $r>0$, the function $X\mapsto n_{r}(X)$ is
    upper semicontinuous on compact sets in the Gromov-Hausdorff
    topology.
    
    \itm{ii} For any $r>0$, $R>0$, the function $(X,x) \mapsto
    n_{r}(\B_{X}(x,R))$ is upper semicontinuous on proper spaces in
    the pointed Gromov-Hausdorff topology.
    
    \item{$(iii)$} For any $r>0$, the function $X\mapsto
    \overline{n}_{r}(X)$ is lower semicontinuous on compact sets in
    the Gromov-Hausdorff topology.
    
    \itm{iv} For any $r>0$, $R>0$, the function $(X,x) \mapsto
    \ov{n}_{r}(\ov{B_{X}(x,R)})$ is lower semicontinuous on proper
    spaces in the pointed Gromov-Hausdorff topology.
\end{prop}

\begin{proof}
    $(i)$.  Since $n_{r}$ is integer valued, the statement is
    equivalent to: $\forall K$ compact $\exists\d>0$ s.t.
     $d_{GH}(J,K)<\d$, imply $n_{r}(J)\leq n_{r}(K)$. 
    Then, if $\cup_{j=1}^{n_{r}(K)}B(x_{j},r)$ is a minimal
    covering for $K$ with open balls of radius $r$, and we set
    $$
    R=\max_{x\in K}\min_{j=1,\dots n_{r}(K)}d_{K}(x,x_{j}),
    $$ 
    then $\d=r-R>0$ and $n_{\r}(K)=n_{r}(K)$ for any $r\geq\r>R$. 
    Therefore $d_{H}(J,K)<\d/2$ implies $J$ and $K$ may be embedded in
    a metric space $Z$ where $J\subset B_{Z}(K,\d/2)$, $K\subset
    B_{Z}(J,\d/2)$, hence we may find points $y_{1},\dots,y_{n}\in J$
    with $d_{Z}(x_{i},y_{i})<\d/2$.  Finally
    $$
    \cup_{j=1}^{n_{r}(K)}B(y_{j},r)\supset 
    \cup_{j=1}^{n_{r}(K)}B(x_{j},r-\d/2)\supset B_{Z}(K,\d/2),
    $$
    namely $n_{r}(J) \leq n_{r}(K)$.
    
    $(ii)$.  Assume  $(X_{n},x_{n})\convPGH(X,x)$.  Then, by 
    Proposition \ref{Prop:Gromovbook} for any
    given $R$, there exists $\eps_{n}\to0$ such that
    $\ov{B}_{X_{n}}(x_{n},R+\eps_{n})\convGH \ov{B}_{X}(x,R)$.  Eventually, 
    by $(i)$, 
    $$
    n_{r}(\ov{B}_{X_{n}}(x_{n},R))\leq 
    n_{r}(\ov{B}_{X_{n}}(x_{n},R+\eps_{n}))
    \leq n_{r}(\ov{B}_{X}(x,R)).  
    $$

    $(iii)$.  We have to show that, for any $p\in\bn$, $\cb := \{ X
    \text{ proper metric space }: \ov{n}_{r}(X)\leq p\}$ is closed. 
    Let $\{X_{n}\}\subset\cb$, $X_{n}\convGH X$, and possibly passing 
    to a subsequence, we may assume that
    $\ov{n}_{r}(X_{n}) = q\leq p$, all $n\in\bn$. 
    According to \cite{BHbook} we may describe $X$ as an ultralimit,
    namely, given any free ultrafilter $\cu$ on $\bn$, we may set
    $d_{\cu}(\{x_{n}\},\{y_{n}\})=\lim_{\cu} d_{X_{n}}(x_{n},y_{n})$
    where $x_{n},y_{n}\in X_{n}$, and $X$ is isometric to the space of
    equivalence classes $x_{\cu}$ of $\{x_{n}\}$ obtained by identifying
    points with zero distance.  Now let $x^{j}_{n}$, $j=1,\dots,q$
    be the centers of balls of radius $r$ covering $X_{n}$, and set 
    $x^{j}_{\cu} := [x^{j}_{n}]$, $j=1,\ldots,q$.  Given any
    $x_{\cu}=[x_{n}]\in X$, setting $N_{j}:= \{n\in\bn: 
    d_{X_{n}}(x_{n},x^{j}_{n})\leq r\}$, there is 
    $j_{0}\in\{1,\ldots,q\}$ such that $N_{j_{0}}\in\cu$, so that 
    $$
    d_{\cu}(x^{j_{0}}_{\cu},x_{\cu}) = \lim_{\cu}d_{X_{n}}
    (x^{j_{0}}_{n},x_{n}) \leq r,
    $$
    hence $\ov{n}_{r}(X)\leq q$, that is $X\in\cb$.
        
    $(iv)$.  Assume $(X_{n},x_{n})\convPGH(X,x)$.  Then, for any given
    $R$, and possibly passing to a subsequence, 
    $\ov{B_{X_{n}}(x_{n},R)}\convGH B$ with $\ov{B_{X}(x,R)}\subset 
    B\subset \ov{B}_{X}(x,R)$ (cf. Proposition \ref{Prop:conseq}). Therefore, 
    by $(iii)$,
    $$
    \ov{n}_{r}\left(\ov{B_{X}(x,R)}\right)\leq \ov{n}_{r}(B)\leq
    \liminf \ov{n}_{r}\left(\ov{B_{X_{n}}(x_{n},R)}\right).
    $$
\end{proof}

 \begin{thm}\label{dim-formula}
     Let $(X,d)$ be a metric space, and let $x\in X$ be such that the 
     sufficient condition (\ref{cptcond}) is satisfied.  The following
     formulas hold:
     \begin{align*}
	\supc_{X}(x)&=\limsup_{r\to 0}\sup_{T\in\ct_{x} X}\frac{\log
	n_{r}(\B_{T}(x,1))}{\log1/r}, \\
	\subc_{X}(x)&=\liminf_{r\to 0}\inf_{T\in\ct_{x} X}\frac{\log
	n_{r}(\B_{T}(x,1))}{\log1/r}.
    \end{align*}
 \end{thm}
 \begin{proof}
     Let us denote by $tX$ the metric space $(X,td)$.
     Fix $\l>0$ and choose $r_{n}\to 0$ such that 
     $\limsup_{r\to0}n_{\l r}(\B_{X}(r))=
     \limsup_{n}n_{\l r_{n}}(\B_{X}(r_{n}))$ and $r_{n}^{-1}X$ is 
     converging in the pointed Gromov-Hausdorff topology, say to a 
     tangent set $T$. By the proposition above, 
     \begin{align*}
	 n_{\l}(B_{T}(1))&\geq\limsup_{n}
	 n_{\l }(\B_{r_{n}^{-1}X}(1))=
	 \limsup_{n}n_{\l r_{n}}(\B_{X}(r_{n}))\\
	 &=\limsup_{r\to0}n_{\l r}(\B_{X}(r)).
     \end{align*}
     Taking the $\limsup_{\l\to0}\sup_{T\in\ct_{x}(X)}$ we get
     $$
     \supc_{X}(x)\leq\limsup_{r\to 0}\sup_{T\in\ct_{x} X}\frac{\log
     n_{r}(\B_{T}(x,1))}{\log1/r}.
     $$
     
     Conversely, for any $T\in \ct_{x}(X)$, with $r_{n}^{-1}X\convPGH 
     T$, we get
     $\ov{n}_{\l}(\ov{B_{T}(1)})\leq \liminf_{n}\ov{n}_{\l 
     r_{n}}(\ov{B_{X}(r_{n})})\leq \limsup_{r}\ov{n}_{\l 
     r}(\ov{B_{X}(r)})$.
     
     Taking the $\limsup_{\l\to0}\sup_{T\in\ct_{x}(X)}$ we get
     $$
     \supc_{X}(x)\geq\limsup_{r\to 0}\sup_{T\in\ct_{x} X}\frac{\log
     \ov{n}_{r}(\ov{B_{T}(x,1)})}{\log1/r}.
     $$
     The thesis easily follows.
 \end{proof}
 
 \subsection{Further properties of tangential dimensions}

 Tangential dimensions are invariant under bi-Lipschitz maps.
 
 \begin{prop}
	 Let $X,\ Y$ be metric spaces, $f:X \to Y$ be a bi-Lipschitz map
	 $i.e.$ there is $L>0$ such that $L^{-1} d_{X}(x,x') \leq
	 d_{Y}(f(x),f(x')) \leq L d_{X}(x,x')$, for $x,x'\in X$.  Then
	 $\underline{\d}_{X}(x) = \underline{\d}_{Y}(f(x))$ and 
	 $\overline{\d}_{X}(x) = \overline{\d}_{Y}(f(x))$, for all $x\in X$.
 \end{prop}
 \begin{proof}
 	Observe that, for any $x\in X$, $y\in Y$, $r>0$, we have 
	\begin{align}\label{eqn:Lip}
		B(f(x),r/L) & \subset f(B(x,r)) \subset B(f(x),rL) \\
		B(f^{-1}(y),r/L) & \subset f^{-1}(B(y,r)) \subset
		B(f^{-1}(y),rL)
	\end{align}
	so that
	$$
	f(B(x,R)) \subset B(f(x),RL) \subset \bigcup_{i=1}^{n(r/L,B(f(x),RL))} 
	B(y_{i},r/L)
	$$
	and
	$$
	B(x,R) \subset \bigcup_{i=1}^{n(r/L,B(f(x),RL))}
	f^{-1}(B(y_{i},r/L)) \subset \bigcup_{i=1}^{n(r/L,B(f(x),RL))} 
	B(f^{-1}(y_{i}),r)
	$$
	from which it follows $n(r,B(x,R)) \leq n(r/L,B(f(x),RL))$. Exchanging 
	the roles of $f$ and $f^{-1}$, we obtain $n(r,B(f(x),R)) \leq 
	n(r/L,B(x,RL))$, so that
	$$
	n(rL,B(f(x),\frac{R}{L})) \leq n(r,B(x,R)) \leq 
	n(r/L,B(f(x),RL)).
	$$
	Therefore, taking $\limsup_{R\to 0}$, then $\limsup_{\l\to0}$, and 
	doing some algebra, we get
	\begin{align*}
	    \limsup_{\l\to0} \limsup_{R\to 0} \frac{\log n(\l
	    R,B(f(x),R))}{\log L^{2}/\l} & \leq \limsup_{\l\to0} \limsup_{R\to
	    0} \frac{\log n(\l R,B(x,R))}{\log 1/\l} \\
	    & \leq \limsup_{\l\to0} \limsup_{R\to 0} \frac{\log
	    n(\l R,B(f(x),R))}{\log 1/(L^{2}\l)}
	\end{align*}
	which means $\overline{\d}_{X}(x) = \overline{\d}_{Y}(f(x))$. The 
	other equality is proved in the same manner.
 \end{proof}

 The following proposition shows that the functions 
 $\underline{\d}_{X}$ and $\overline{\d}_{X}$ satisfy properties which 
 are characteristic of a dimension function. Denote by $B_{Y}(x,r):= 
 Y\cap B_{X}(x,r)$, if $Y\subset X$.
 
 \begin{prop}
     \itm{i} Let $Y\subset X$ and $x\in Y$. Then 
     $\underline{\d}_{Y}(x) \leq \underline{\d}_{X}(x)$, and 
     $\overline{\d}_{Y}(x) \leq \overline{\d}_{X}(x)$. Equality holds 
     if there is $R_{0}>0$ such that $B_{X}(x,R_{0}) \subset Y$.
     \itm{ii} Let $X_{1},\ X_{2} \subset X$ and $x\in X_{1}\cap 
     X_{2}$. Then
     \begin{align*}
	 \underline{\d}_{X_{1}\cup X_{2}}(x) & \geq \max \{
	 \underline{\d}_{X_{1}}(x), \underline{\d}_{X_{2}}(x) \} \\
	 \overline{\d}_{X_{1}\cup X_{2}}(x) & = \max \{
	 \overline{\d}_{X_{1}}(x), \overline{\d}_{X_{2}}(x) \}.
     \end{align*}
     \itm{iii} Let $X,\ Y$ be metric spaces, $x\in X$, $y\in Y$.  Then
     $\underline{\d}_{X\times Y}((x,y)) \geq \underline{\d}_{X}(x) +
     \underline{\d}_{Y}(y)$, and $\overline{\d}_{X\times Y}((x,y)) \leq
     \overline{\d}_{X}(x) + \overline{\d}_{Y}(y)$.
 \end{prop}
 \begin{proof}
     $(i)$ As $B_{Y}(x,R) \subset B_{X}(x,R)$, we get
     $n_{r}(B_{Y}(x,R)) \leq n_{r}(B_{X}(x,R))$, and analogously for
     $\n_{x}$, and the claim follows.  The second statement is
     obvious.
     
     $(ii)$ The inequalities $\geq$ follow from $(i)$.  It remains to
     prove $\overline{\d}_{X_{1}\cup X_{2}}(x) \leq \max \{
     \overline{\d}_{X_{1}}(x), \overline{\d}_{X_{2}}(x) \}$, and we 
     can assume $a:= \overline{\d}_{X_{1}}(x)<\infty$ and $b:= 
     \overline{\d}_{X_{2}}(x) <\infty$, otherwise there is nothing to 
     prove. Assume for definiteness that $a\leq b$. Then, for any 
     $\eps>0$, there is $\l_{0}>0$ such that, for any 
     $\l\in(0,\l_{0})$, there exists $r_{0} = r_{0}(\eps,\l)$ such 
     that, for any $r\in(0,r_{0})$ we get
     \begin{align*}
	 n_{\l r}(B_{X_{1}}(x,r)) \leq \frac{1}{\l^{a+\eps}} \\
	 n_{\l r}(B_{X_{2}}(x,r)) \leq \frac{1}{\l^{b+\eps}}. 
     \end{align*}
     As $B_{X_{1}\cup X_{2}}(x,R) \subset B_{X_{1}}(x,R) \cup
     B_{X_{2}}(x,R)$, we get 
     \begin{align*}
	 n_{\l r}(B_{X_{1}\cup X_{2}}(x,r)) & \leq n_{\l
	 r}(B_{X_{1}}(x,r)) + n_{\l r}(B_{X_{2}}(x,r)) \\
	  & \leq \frac{1}{\l^{a+\eps}} + \frac{1}{\l^{b+\eps}} \\
	  & = \frac{1}{\l^{b+\eps}} (1 + \l^{b-a}).
     \end{align*}
     Therefore
     $$
     \frac{\log n_{\l r}(B_{X_{1}\cup X_{2}}(x,r))}{\log 1/\l} \leq 
     b + \eps + \frac{\log(1 + \l^{b-a})}{\log 1/\l},
     $$
     so that $\overline{\d}_{X_{1}\cup X_{2}}(x) \leq b+\eps$, and 
     the thesis follows by the arbitrariness of $\eps$.
     
     $(iii)$ Endow $X\times Y$ with the metric
     \begin{equation}\label{eqn:prod.metric}
	 d((x_{1},y_{1}),(x_{2},y_{2})) := \max \{ d_{X}(x_{1},x_{2}),
	 d_{Y}(y_{1},y_{2}) \}
     \end{equation}
     which is by-Lipschitz equivalent to the
     product metric.  Then 
     \begin{equation}\label{eqn:prod.ball}
	 B_{X\times Y}((x,y),R) = B_{X}(x,R) \times B_{Y}(y,R).
     \end{equation}
     Therefore $\n_{r}(B_{X\times Y}((x,y),R)) \geq \n_{r}(B_{X}(x,R))
     \n_{r}(B_{Y}(y,R))$, and
     \begin{align*}
	 \liminf_{\l\to0} \liminf_{r\to0} \frac{\log \n_{\l
	 r}(B_{X\times Y}((x,y),r))}{\log 1/\l} & \geq
	 \liminf_{\l\to0} \liminf_{r\to0} \frac{\log \n_{\l
	 r}(B_{X}(x,r))}{\log 1/\l} \\
	 & + \liminf_{\l\to0} \liminf_{r\to0} \frac{\log \n_{\l
	 r}(B_{Y}(y,r))}{\log 1/\l}.
     \end{align*}
     Moreover $n_{r}(B_{X\times Y}((x,y),R)) \leq
     n_{r}(B_{X}(x,R)) n_{r}(B_{Y}(y,R))$, and
     \begin{align*}
	 \limsup_{\l\to0} \limsup_{r\to0} \frac{\log n_{\l
	 r}(B_{X\times Y}((x,y),r))}{\log 1/\l} & \leq
	 \limsup_{\l\to0} \limsup_{r\to0} \frac{\log n_{\l
	 r}(B_{X}(x,r))}{\log 1/\l} \\
	 & + \limsup_{\l\to0} \limsup_{r\to0} \frac{\log n_{\l
	 r}(B_{Y}(y,r))}{\log 1/\l}.
     \end{align*}
 \end{proof}

 \section{Local dimensions of tangent sets}\label{section:g&dg}

 \subsection{A different formula for tangential
 dimensions}\label{subsect:comparison}

 There is another notion of dimension naturally associated with the
 tangent cone.  One may indeed take the infimum, resp.  supremum, of
 the lower, resp.  upper, box dimension of the tangent balls at a
 given point.  We give below a sufficient condition for them to coincide
 with the tangential dimensions defined above.  However, this equality
 does not hold in general, as shown in subsection \ref{ex:counterexmp}.

 Let $X$ be a metric space, $x\in X$.  We shall consider
 the following.

 \begin{Assump}\label{Ass:bdd}
    There exist constants $c\geq1$, $a\in(0,1]$ such that, for any
    $r\leq a$, $\l,\m\leq1$, $y,z\in B_{X}(x,r)$,
    \begin{equation}\label{eq:ass}
	n(\l\m r, B_{X}(y,\l r))\leq c n(\l\m r, B_{X}(z,\l r)).
    \end{equation}
 \end{Assump}
 
 Let us observe that the previous inequality is trivially satisfied 
 when $\m\geq1$. 

 Let us recall some notions of dimension.

 The lower and upper box dimensions of $X$ are 
 \begin{align*}
    \subd(X) & = \lim_{R \to \infty} \liminf_{r\to0} \frac{\log
    n(r,B_{X}(x,R))}{\log 1/r}, \\
    \supd(X) & = \lim_{R \to \infty} \limsup_{r\to0} \frac{\log
    n(r,B_{X}(x,R))}{\log 1/r},
 \end{align*}
 while the (lower and upper) local (box) dimensions of $X$ at a point $x$
 are defined as
 \begin{align*}
    \subd_{X}(x)&=\lim_{R\to0}\subd(B_{X}(x,R)) = \lim_{R \to
    0} \liminf_{r\to0} \frac{\log n(r,B_{X}(x,R))}{\log 1/r},\\
    \supd_{X}(x)&=\lim_{R\to0}\supd(B_{X}(x,R)) = \lim_{R \to
    0} \limsup_{r\to0} \frac{\log n(r,B_{X}(x,R))}{\log 1/r}.
 \end{align*}
 
 \begin{rmk}\label{rem:equiv-form}
     $(i)$ For the box dimensions to be non-trivial, (the completion 
     of) $X$ has to be proper; for the local box dimensions at $x$ 
     to be non-trivial, $x$ needs to have a compact (totally bounded) 
     neighborhood.
     
     $(ii)$ We obtain the same definition if we replace $n$ 
     with $\n$ or with $\ov{n}$, and/or $B_{X}(x,r)$ with 
     $\B_{X}(x,r)$ or with $\ov{B_{X}(x,r)}$. The proof is the 
     same as that of Proposition \ref{equiv-form}.
 \end{rmk}
 
 Now we set
 \begin{equation}\label{g(t,h)}
     g(t,h)=\log n(e^{-(t+h)},B_{X}(x,e^{-t})).
 \end{equation}
 Clearly $h\mapsto g(t,h)$ is non-decreasing for any $t$.  The
 tangential dimensions can be rewritten as
 \begin{align*}
    \subc_{X}(x) 
    &= \liminf_{h\to+\infty} \liminf_{t\to+\infty}\frac{g(t,h)}{h},\\
    \supc_{X}(x) 
    &= \limsup_{h\to+\infty} \limsup_{t\to+\infty}\frac{g(t,h)}{h}.
 \end{align*}
 
 The local dimensions can be rewritten as
 \begin{align*}
    \subd_{X}(x) 
    &= \lim_{t\to+\infty} \liminf_{h\to+\infty}\frac{g(t,h)}{h},\\
    \supd_{X}(x) 
    &= \lim_{t\to+\infty} \limsup_{h\to+\infty}\frac{g(t,h)}{h}.
 \end{align*}

We define the {\it coboundary} of $g$ as the three-variable function
$$
dg(t,h,k)=g(t, h+k) - g(t+h,k)-g(t,h),
$$
and note that $g$ is a {\it cocycle}, namely $dg=0$, if and only if 
$g(t,h)=g(0,t+h)-g(0,t)$, namely if it is a coboundary where, given 
$t\to f(t)$, we set $df(t,h)=f(t+h)-f(t)$.

We shall show that our assumption implies a bound on $dg$. 

\begin{lem}\label{number-ineq}
    The following inequality holds:
    \begin{equation*}
	n(\l\m r,B_{X}(x,r))
	\leq n(\l r,B_{X}(x,r)) 
	\sup_{y\in B_{X}(x,r)}n(\l\m r,B_{X}(y,\l r)).
    \end{equation*}
\end{lem}

\begin{proof}
   Let us note that we may realize a covering of $B_{X}(x,r)$
   with balls of radius $\l\m r$ as follows: first choose an
   optimal covering of $B_{X}(x,r)$ with balls of radius
   $\l r$, and then cover any covering ball optimally with balls
   of radius $\l\m r$.  The thesis follows.    
\end{proof}

Let us recall that the function $\n(r,B_{X}(x,R))$, denotes the
maximum number of disjoint open balls of $X$ of radius $r$ centered in
the open ball of center $x$ and radius $R$ of $X$.

\begin{lem}\label{ni-ineq}
    The following inequality holds:
    \begin{equation*}
	\n(\l\m r,B_{X}(x,r))\geq 
	\n(\l r,B_{X}(x,r)) 
	\inf_{y\in B_{X}(x,r)}\n(\l\m r,B_{X}(y,\l r)).
    \end{equation*}
\end{lem}

\begin{proof}
    Indeed we may find disjoint open balls of $X$ of radius $\l\m r$
    centered in $B_{X}(x,r)$ as follows: first find a maximal set of
    disjoint open balls of $X$ of radius $\l r$ centered in
    $B_{X}(x,r)$, and then, for any such ball, find a maximal set of
    disjoint open balls of $X$ of radius $\l\m r$ centered in it.  This
    implies the thesis.
\end{proof}

 \begin{prop}
    Assumption \ref{Ass:bdd} and condition $(\ref{cptcond})$ for $(X,x)$
    imply that $dg$ is bounded, for $t>t_{0}$, $h,k>0$.
 \end{prop}
 \begin{proof}
     Let us observe that it is enough to find a bound for $h,k$
     sufficiently large.  By the assumption and Lemma
     \ref{number-ineq}, for $r\leq a$,
     \begin{align*}
	 n(\l\m r,B_{X}(x,r))
	 &\leq n(\l r,B_{X}(x,r)) \sup_{y\in
	 B_{X}(x,r)}n(\l\m r,B_{X}(y,\l r))\\
	 &\leq c n(\l r,B_{X}(x,r))
	 n(\l\m r,B_{X}(x,\l r)).
     \end{align*}
     Therefore, if we set $r=e^{-t}$, $\l=e^{-h}$,
     $\m=e^{-k}$, we get, for $t\geq \log 1/a$,
     \begin{equation}\label{upperbound}
	 dg(t,h,k)=g(t, h+k) - g(t+h,k)-g(t,h) \leq\log c.
     \end{equation}
     Let us now find a bound from below.  By Lemma \ref{ni-ineq}, the
     inequalities (\ref{n-ni}), and assumption \ref{eq:ass}, we get,
      for $r\leq a$,
    \begin{align*}
	n(\l\m r,B_{X}(x,r))
	&\geq\n(\l\m r,B_{X}(x,r))\\
	& \geq \n(\l r ,B_{X}(x,r)) \inf_{y\in
	B_{X}(x,r)}\n(\l\m r,B_{X}(y,\l r))\\
	&\geq n(2\l r,B_{X}(x,r))
	\inf_{y\in B_{X}(x,r)}n(2\l\m r,B_{X}(y,\l r))\\
	&\geq \frac{1}{c}n(2\l r,B_{X}(x,r))
	n(2\l\m r,B_{X}(x,\l r)).
    \end{align*}
    As a consequence, for $t\geq \log 1/a$,
    $$
    g(t, h+k)-g(t,h-\log2)-g(t+h,k-\log2)\geq-\log c,
    $$
    which implies
    $$
    dg(t,h,k)\geq-\log c -( g(t,h)-g(t,h-\log2))-
    (g(t+h,k)-g(t+h,k-\log2)).
    $$
    The result follows if we show that $g(t,h+\log2)-g(t,h)$ is
    bounded from above. Indeed, by the upper bound (\ref{upperbound}),  
    \begin{align*}
	g(t,h+\log2)-g(t,h)&=dg(t,h,\log2)+g(t+h,\log2)\\
	&\leq \log{c}+g(t+h,\log2).
    \end{align*}
    Since $L=\limsup_{t\to +\infty}g(t,\log 2)<+\infty$, there exists 
    $t_{0}>\log 1/a$ such that, for $t>t_{0}$, $g(t,\log 2)\leq 2 L$.
\end{proof}

\begin{prop}\label{g-prop}
    Let us assume property \ref{Ass:bdd}.  
    \item{$(i)$} If condition (\ref{cptcond}) holds, the
    $\liminf_{\l\to 0}$, resp.  $\limsup_{\l\to 0}$, in the definition
    of $\subc$, resp.  $\supc$, are indeed limits:
    \begin{align*}
	\underline{\d}_{X}(x) & = \lim_{\l \to 0} \liminf_{r \to
	0} \frac{\log n(\l r,B(x,r))}{\log 1/\l},
	\\
	\overline{\d}_{X}(x) & = \lim_{\l \to 0} \limsup_{r \to
	0} \frac{\log n(\l r,B(x,r))}{\log 1/\l}.
    \end{align*}
    \item{$(ii)$} If condition (\ref{cptcond}) holds, the
    following inequalities hold:
    $$
    \subc_{X}(x)\leq\subd_{X}(x)\leq\supd_{X}(x)\leq\supc_{X}(x).
    $$
    \item{$(iii)$} Condition (\ref{cptcond}) is equivalent 
    to the finiteness of $\overline{\d}_{X}(x)$.
\end{prop}
\begin{proof}
    All the statements follow directly by Proposition \ref{lim&ineq}.
\end{proof}

\begin{lem}\label{dimtang}
    Let $\l_{n}\to0$ be a sequence such that 
    $(\frac{1}{\l_{n}}X,x)\convPGH (T,x)$. Then,
    \begin{align}
	\subd(\B_{T}(x,1))&=\liminf_{h\to\infty}\liminf_{n}\frac{g(t_{n},h)}{h}
	=\liminf_{h\to\infty}\limsup_{n}\frac{g(t_{n},h)}{h},\label{dimtg1}\\
	\supd(\B_{T}(x,1))&=\limsup_{h\to\infty}\liminf_{n}\frac{g(t_{n},h)}{h}
	=\limsup_{h\to\infty}\limsup_{n}\frac{g(t_{n},h)}{h},\label{dimtg2}
    \end{align}
    where we posed $t_{n}=-\log\l_{n}$.
\end{lem}
\begin{proof}
    In the following we shall omit the reference to the point $x$.  By
    definition, setting $h=\log1/r$,
    \begin{align*}
	\limsup_{n}\frac{g(t_{n},h)}{h}
	&=\limsup_{n}\frac{\log n(\l_{n}r,\B_{X}(\l_{n}))}{\log1/r}\\
	&=\limsup_{n}\frac{\log n(r,\B_{1/\l_{n}X}(1))}{\log1/r}\\
	&\leq \frac{\log n(r,\B_{T}(1))}{\log1/r},
    \end{align*}
    where we used the upper semicontinuity in Proposition
    \ref{semicontBis} $(ii)$.  Analogously,
    \begin{align*}
	\liminf_{n}\frac{\log \overline{n}(\l_{n}r,\ov{
	B_{X}(\l_{n}) })}{\log1/r} &=\liminf_{n}\frac{\log
	\overline{n}(r,\ov{ B_{1/\l_{n}X}(1) })}{\log1/r}\\
	&\geq \frac{\log \overline{n}(r,\ov{ B_{T}(1) })}{\log1/r},
    \end{align*}
    namely
    \begin{align*}
	\frac{\log \overline{n}(r,B_{T}(1))}{\log1/r}
	& 
	\leq \frac{\log \overline{n}(r,\ov{ B_{T}(1) })}{\log1/r} \\
	& \leq
	\liminf_{n}\frac{\log \overline{n}(\l_{n}r,\ov{B_{X}(\l_{n}) } )}
	{\log1/r}
	\\
	& \leq
	\limsup_{n}\frac{\log n(\l_{n}r,\B_{X}(\l_{n}) )}{\log1/r} \\
	& \leq
	\frac{\log n(r,\B_{T}(1))}{\log1/r}.
    \end{align*}
    Recalling Proposition \ref{equiv-form} and Remark
    \ref{rem:equiv-form} $(ii)$, and taking the $\liminf$ for $r\to0$ we get
    the equalities (\ref{dimtg1}), taking the $\limsup$ for $r\to0$ we
    get the equalities (\ref{dimtg2}).
\end{proof}

\begin{thm}\label{newformula}
    Under the Assumption \ref{Ass:bdd} and condition (\ref{cptcond})
    \begin{align}
	\subc_{X}(x) &= \inf_{T\in\ct_{x} X}\subd(T)=
	\inf_{T\in\ct_{x} X}\supd(T),\label{eq.5}\\
	\supc_{X}(x) &= \sup_{T\in\ct_{x} X}\subd(T)=
	\sup_{T\in\ct_{x} X}\supd(T).\label{eq.6}
    \end{align}
\end{thm}

\begin{proof}
    We only prove (\ref{eq.6}), the proof of (\ref{eq.5}) being
    analogous.  Let us observe that the property satisfied by the
    sequence $\overline{t}_{n}\to\infty$ described in Proposition
    \ref{minimizing} remains valid for any subsequence.  We may
    therefore assume that $\overline{t}_{n}$ produces a tangent set,
    namely $e^{\overline{t}_{n}}X$ converges to a tangent set $T$ in
    the pointed Gromov-Hausdorff topology.  Then, by Lemma
    \ref{dimtang} and Proposition \ref{minimizing}, for any $\k$ there
    exists a tangent set $T$ such that $\supc_{X}(x) -
    \frac{2S}{\k}\leq \supd(\B_{T}(1)) \leq \supc_{X}(x)$, hence
    $$
    \supc_{X}(x) = \sup_{T\in\ct_{x} X}\supd(\B_{T}(1)).
    $$
    Since $\ct_{x} X$ is globally dilation invariant, and the box 
    dimensions are dilation invariant, for any tangent set $T$ and 
    any $r>0$ there exists a tangent set $S$ for which 
    $$
    \supd(\B_{T}(r)) = \supd(r\B_{S}(1))) = \supd(\B_{S}(1)),
    $$
    namely
    $$
    \sup_{T\in\ct_{x} X}\supd(\B_{T}(1)) = \sup_{\substack{T\in\ct_{x} 
    X\\ r>0}}\supd(\B_{T}(r)) = \sup_{T\in\ct_{x} X}\supd(T).
    $$
\end{proof}

\subsection{A counterexample}\label{ex:counterexmp}
Now we show that the equality shown above under hypothesis
\ref{Ass:bdd} does not hold in general.  In the example below we
construct a subset of $\br^{3}$ for which any tangent set at a given
point is zero-dimensional, but $\supc$ is positive.

First set $a^{k}_{n} := \e{-\left( \frac{(n+k)(n+k+1)}{2}+k
\right)^{2}}$, for $k,\,n\in\bn$.  

Let now $S^{2}:= \{ x\in\br^{3} : \|x\|=1 \}$, and choose, for any
$k\in\bn$, $S_{k}\subset S^{2}$ such that 
\begin{itemize}
    \item The diameter of $S_{k}$ is $1/k^{2}$,
    \item $\# S_{k}=k^{2}$,
    \item $d(v,w)\geq\frac{1}{ k^{3} }$, $v,w\in S_{k}$, $v\neq w$,
    \item $\min\{d(v,w):v\in S_{k}, w\in S_{h}\} \geq \frac{1}{ k^{3} }$, 
    $h\geq k$,
    \item $\lim_{k\to\infty} S_{k} =: S_{\infty} =\{v_{\infty} \}\subset 
    S^{2}$ in the Hausdorff topology. 
\end{itemize}

Set, for any $k\in\bn$, $A_{k} := \{ a^{k}_{n}v : v\in S_{k}, n\in\bn
\}$, and $F := \overline{\cup_{k=1}^{\infty} A_{k}} \subset \br^{3}$. 

\begin{lem}\label{Finitetangent}
    The tangent cone of $F$ at $0$ consists, up to dilations, of the 
    set $\{0\}$ and of the sets 
    $S_{k}\cup\{0\}$, with $k\in\bn\cup\{\infty\}$.
\end{lem}
\begin{proof}
    If $(\l_{n}F,0)\convPGH (T,0)$, then the tangent set $T$ does not
    consist of the sole $\{0\}$ if and only if, for suitable sequences
    $n(p)$, $k(p)\in\bn$, $v(p)\in S_{k(p)}$, $\l_{p}a^{k(p)}_{n(p)}
    v(p)$ converges, when $p\to\infty$, and $\l_{p}a^{k(p)}_{n(p)}\to
    c\in (0,\infty)$.  \\
    Assume $\{k(p)\}$ is bounded.  Since $\l_{p} a^{k(p)}_{n(p)} v(p)$
    converges, then $k(p)$ has to be eventually equal to some $k_{0}$,
    namely we may replace $\l_{p}$ with a subsequence of $c\left(
    a^{k_{0}}_{n} \right)^{-1}$.  This implies that $T\supseteq
    cS_{k}$.  \\
    Let us observe that two infinitesimal subsequences $c_{n}$,
    $c'_{n}$ contained in $\{a^{k}_{n}:k,n\in\bn\}$ such that
    $\frac{c_{n}}{c'_{n}}\to \c\ne 0$ eventually coincide.
    From this it is not difficult to derive that all limit points in
    $(\l_{n}F,0)$ belong to $cS_{k}$.
    \\
    If $\{k(p)\}$ is not bounded, it has to diverge, namely
    $T\supseteq cS_{\infty}$.  Reasoning as before, one gets $T=cS_{\infty}$.
\end{proof}

 \begin{prop} \label{prop:counterexmp3}
     Let $F$ be as above.  Then
     $$\supc_{F}(0) > \sup_{T\in\ct_{0}F} \supd(T) = 0.$$
 \end{prop}
 \begin{proof}
     By Lemma \ref{Finitetangent} we get $\sup_{T\in\ct_{0}F} \supd(T) = 0$.
     Now let $k\in\bn$, and let $\{\l_{n}\}\subset(0,\infty)$ be an
     increasing diverging sequence s.t. $X:=
     \lim_{n\to\infty} \l_{n} F$ exists and $\B_{X}(0,1)$
     consists of $k+1$ points, all belonging to $\{tv : v\in S_{k}, 
     t\geq0\}$. As $n_{1/k^{2}}(\B_{X}(0,1)) = k+1$, we obtain
     $$
     \frac{\log n_{1/k^{2}}(\B_{X}(0,1))}{\log k^{2}} \geq \frac12,
     $$
     so that 
     $$
     \supc_{F}(0) = \limsup_{r\to0} \sup_{T\in\ct_{0}F} \frac{\log 
     n_{r}(\B_{T}(0,1))}{\log 1/r} \geq \frac12.
     $$
 \end{proof}

 \section{Appendix}

 Here we collect some results on the two-variable functions $g(t,h)$. 
 Throughout this section we assume that $g$ is non-decreasing in the 
 $h$ variable and that, for a suitable constant $t_{0}$,
 $$
 S=\sup_{
 \begin{smallmatrix}
     t>t_{0}\\
     h,g>0
 \end{smallmatrix}}|dg(t,h,k)|<\infty,
 $$
 where $dg(t,h,k)=g(t, h+k) - g(t+h,k)-g(t,h)$.

\begin{lem}
    Given $t>t_{0}$, $h_{1},\dots h_{n}>0$, we have
    \begin{equation}\label{quasicocycle}
	\left|g(t,\sum_{i=1}^{n}h_{i})-
	\sum_{k=1}^{n}g(t+\sum_{i=1}^{k-1}h_{i},h_{k})\right|
	\leq (n-1)S.
    \end{equation}
\end{lem}
\begin{proof}
    A straightforward computation gives
    \begin{equation}
	g(t,\sum_{i=1}^{n}h_{i})=
	\sum_{k=1}^{n}g(t+\sum_{i=1}^{k-1}h_{i},h_{k})+
	\sum_{k=1}^{n-1}dg(t+\sum_{i=1}^{k-1}h_{i},h_{k},
	\sum_{i=k+1}^{n}h_{i}).
    \end{equation}
    The thesis follows.
\end{proof}

\begin{prop}\label{lim&ineq}
    \item{$(i)$} The quantities
    $$
    \limsup_{t\to\infty}\frac{g(t,h)}{h},\quad
    \liminf_{t\to\infty}\frac{g(t,h)}{h},
    $$
    have a limit when $h\to\infty$.
    \item{$(ii)$} The following inequalities hold:
    \begin{align*}
	\lim_{t\to\infty}\limsup_{h\to\infty}\frac{g(t,h)}{h}
	&\leq\lim_{h\to\infty}\limsup_{t\to\infty}\frac{g(t,h)}{h},\\
	\lim_{t\to\infty}\liminf_{h\to\infty}\frac{g(t,h)}{h}
	&\geq\lim_{h\to\infty}\liminf_{t\to\infty}\frac{g(t,h)}{h}.	
    \end{align*}
    \item{$iii$} The quantity 
    $\lim_{h\to\infty}\limsup_{t\to\infty}\frac{g(t,h)}{h}$ is 
    infinite if and only if the quantity
    $\limsup_{t\to\infty}g(t,h)$ is infinite for one (and in fact for 
    any) $h>0$.
\end{prop}
\begin{proof}
    $(i)$.
    Let us set $\subg(h)=\liminf_{t\to\infty}g(t,h)$. Then, by eq. 
    (\ref{quasicocycle}), we get
    \begin{equation}\label{h-nh}
	\frac{\subg(nh)}{nh}\geq\frac{\subg(h)}{h}-\frac{S}{h}.
    \end{equation}
    Therefore,
    $$
    \frac{\subg(s)}{s}\geq
    \frac{\subg\left(\floor{\frac{s}{r}}r\right)}{s}
    \geq\floor{\frac{s}{r}}\frac{r}{s}
    \left(\frac{\subg(r)}{r}-\frac{S}{r}\right).
    $$
    Taking the $\liminf_{s\to\infty}$, we get
    \begin{equation}\label{quasi-ineq}
	\liminf_{s\to\infty}\frac{\subg(s)}{s}\geq
	\frac{\subg(r)}{r}-\frac{S}{r}.
    \end{equation}
    Then we take the $\limsup_{r\to\infty}$, and obtain
    $$
    \liminf_{s\to\infty}\frac{\subg(s)}{s}\geq
    \limsup_{r\to\infty}\frac{\subg(r)}{r},
    $$
    which proves the existence of 
    $\lim_{h\to\infty}\liminf_{t\to\infty}\frac{g(t,h)}{h}$.
    The existence of the other limit is proved analogously.
    \\
    $(ii)$. Since $g$ is non-decreasing in $h$, for any $\k>0$ we have
    $$
    \liminf_{h\to\infty}\frac{g(t,h)}{h}=
    \liminf_{n\in\bn}\frac{g(t,n\k)}{n\k}.
    $$
    Then, again by eq. (\ref{quasicocycle}), we get, for $t>t_{0}$,
    \begin{equation}
	\frac{g(t,n\k)}{n\k}\geq
	\frac{1}{\k} \left(\frac{1}{n} \sum_{k=1}^{n}g(t+(k-1)\k,\k)-S. 
	\right)
    \end{equation}
    Taking the $\liminf$ on $n\in\bn$ we get
    $$
    \liminf_{n\in\bn} \frac{g(t,n\k)}{n\k}\geq
    \liminf_{t\to\infty}\frac{g(t,\k)}{\k}-\frac{S}{\k},
    $$
    from which 
    $$
    \lim_{t\to\infty}\liminf_{h\to\infty}\frac{g(t,h)}{h}
    \geq\lim_{h\to\infty}\liminf_{t\to\infty}\frac{g(t,h)}{h}
    $$
    follows. The other inequality is proved in the same way.
    \\
    $(iii)$.  Sufficiency is obvious.  Conversely,
    set $\displaystyle{\supg(h)=\limsup_{t\to\infty}g(t,h)}$.  Then, by
    eq.  (\ref{quasicocycle}), and in analogy with (\ref{h-nh}), we get
    \begin{equation}\label{h-nh.sup}
	\frac{\supg(nh)}{nh}\leq\frac{\supg(h)}{h}+\frac{S}{h},
    \end{equation}
    hence, taking the $\lim_{n\to\infty}$,
    $$
    \supg(h) \geq h \lim_{h'\to\infty}
    \limsup_{t\to\infty}\frac{g(t,h')}{h'}-S,
    $$
    from which the thesis follows.
\end{proof}

In the following $\k$ is a given positive number, and we set
$p(t,h)=g(t,h)/h$.

\begin{lem}\label{supVlimsup}
    Let us define
    \begin{align*}
	V_{h}^{d}&=\{t>0:p(t,h)>d\},\\ 
	V^{d}&=\{h\in \k\bn:\sup V_{h}^{d}=+\infty\},\\ 
	V&=\{d\in\br:\sup V^{d}=+\infty\}.
    \end{align*}
    Then,
    $$
    \sup V=\limsup_{h\in \k\bn} \limsup_{t\to+\infty}p(t,h).
    $$
\end{lem}

\begin{proof}
    Let us observe that if $L=\limsup_{x\to\infty} f(x)$, we have
    $$
    L=\sup\{T\in\br:\{x\in\br: f(x)>T\}\text{\ is\ unbounded}\}.
    $$
    Then, setting
    $$
    U^{d}=\{h\in \k\bn:\limsup_{t\to\infty}p(t,h)>d\},\qquad
    U=\{d:\sup U^{d}=+\infty\},
    $$
    we have
    $$
    \limsup_{h\in \k\bn} \limsup_{t\to+\infty}p(t,h)=\sup U,
    $$
    and
    $$
    \limsup_{t\to+\infty}p(t,h)=\sup\{d:\sup V_{h}^{d}=+\infty\},
    $$
    hence
    $$
    \limsup_{t\to+\infty}p(t,h)>d\Rightarrow
    \sup V_{h}^{d}=+\infty\Rightarrow
    \limsup_{t\to+\infty}p(t,h)\geq d,
    $$
    which implies
    $$
    U^{d}\subseteq\{h\in \k\bn:\sup V_{h}^{d}=+\infty\}\subseteq
    \bigcap_{\eps>0}U^{d-\eps}.
    $$
    Finally,
    $$
    \sup U^{d}=+\infty\Rightarrow
    \sup V^{d}=+\infty\Rightarrow
    \sup U^{d-\eps}=+\infty,\ \forall\eps>0,
    $$
    from which the thesis follows.
\end{proof}

\begin{lem}\label{Vtilde}
    Let us define
    \begin{align*}
	\tilde{V}_{h}^{d}&=\{t>0:p(t,j)>d,j\in \k\bn,j\leq h\},\\
	\tilde{V}^{d}&=\{h\in \k\bn:\sup \tilde{V}_{h}^{d}=+\infty\},\\
	\tilde{V}&=\{d\in\br:\sup \tilde{V}^{d}=+\infty\}.
    \end{align*}
    Then $0\leq\sup V-\sup \tilde{V}\leq\frac{2S}{\k}$.
\end{lem}
\begin{proof}
    Since $\tilde{V}_{h}^{d}\subset V_{h}^{d}$, we have $\sup
    \tilde{V}\leq\sup V$.  Let us assume that
    $\sup \tilde{V}<d_{1}<d_{2}<\sup V$, for suitable constants
    $d_{1},d_{2}$.  Now $d_{1}\notin\tilde{V}$, hence there exists
    $\overline{h}\in\k\bn$ such that $\sup
    \tilde{V}_{\overline{h}}^{d_{1}}<+\infty$, namely
    \begin{equation}\label{jt}
	\exists\overline{t}: \forall t>\overline{t}, \exists
	j_{t}\in\k\bn, j_{t}\leq \overline{h}:p(t,j_{t})\leq d_{1}.
    \end{equation}
    Also, $d_{2}\in V$, hence we may find $\tilde{h}\in\k\bn$ such
    that $\sup V_{\tilde{h}}^{d_{2}}=+\infty$ and so large that
    $$
    \tilde{h}>\frac{2d_{1}}{d_{2}-d_{1}} \ov{h}.
    $$
    Therefore we may find $t_{0}>\overline{t}$ such that 
    $p(t_{0},\tilde{h})>d_{2}$. 
    \\
    By equation (\ref{jt}), we can now construct inductively a sequence 
    $j_{i}\in\k\bn$, $j_{i}\leq\overline{h}$, such that, setting
    $$
    t_{k}=t_{0}+\sum_{i=1}^{k}j_{i},
    $$
    we get $p(t_{k},j_{k+1})\leq d_{1}$.
    Since $t_{n}\geq t_{0}+n\k$, there exists $\overline{n}\in\bn$ 
    such that
    $$
    t_{\overline{n}}-\overline{h}\leq t_{\overline{n}-1}\leq t_{0}+\tilde{h}< 
    t_{\overline{n}}.
    $$
    Now, by equation (\ref{quasicocycle}), one gets
    \begin{align*}
	d_{2}&< 
	p(t_{0},\tilde{h})\leq\frac{g(t_{0},\sum_{i=1}^{\overline{n}}j_{i})} 
	{\tilde{h}}\\
	&\leq\frac{1}{\tilde{h}} \sum_{k=1}^{\overline{n}} j_{k}
	p(t_{k-1},j_{k})+\frac{\overline{n}-1}{\tilde{h}} S\\
	&\leq \frac{\sum_{i=1}^{\overline{n}}j_{i}}{\tilde{h}}d_{1}
	+\frac{ S}{\k}\\
	&\leq\left(1+\frac{\overline{h}}{\tilde{h}}\right)d_{1}+\frac{ S}{\k}
	\leq \frac{d_{2}+d_{1}}{2}+\frac{ S}{\k}.
    \end{align*}
    The thesis follows.
\end{proof}

\begin{prop}\label{minimizing}
    For any sequence $t_{n}\to\infty$,
    \begin{equation}\label{supinf1}
	\limsup_{h\to+\infty}
	\limsup_{n\in\bn}\frac{g(t_{n},h)}{h}
	\leq\limsup_{h\to+\infty} \limsup_{t\to+\infty}\frac{g(t,h)}{h}.
    \end{equation}
    Moreover, for any $\k>0$, there exists a sequence
    $\{\overline{t}_{n}\}\to\infty$ for which
    \begin{equation}\label{supinf2}
	\limsup_{h\to+\infty} \limsup_{t\to+\infty}\frac{g(t,h)}{h}
	\leq \liminf_{h\to+\infty}
	\liminf_{n\in\bn}\frac{g(\overline{t}_{n},h)}{h}+\frac{2S}{\k}.
    \end{equation}
\end{prop}
\begin{proof}
    The first inequality is obvious. We shall prove the second.
    \\
    For any given $\k>0$, let $d<\sup\tilde{V}$.  Then $\sup
    \tilde{V}^{d}=+\infty$, i.e. there is $\{h_{n}\}\subset \k\bn$,
    $h_{n}\to\infty$, such that $\sup\tilde{V}^{d}_{h_{n}}=+\infty$. 
    It is not restrictive to assume $h_{n}>n$.  Correspondingly we
    find sequences $t_{nk}\to+\infty$ for $k\to\infty$ such that
    $$
    p(t_{nk},j)>d,\quad j\leq h_{n},\ j\in\k\bn.
    $$
    Again, it is not restrictive to assume $t_{nk}>k$.  Now we make
    explicit the dependence on $d$, setting $h_{np}$ for the sequence
    $h_{n}$ associated to $d=\sup\tilde{V}-1/p$, and $t_{nkp}$ for the
    sequence $t_{nk}$ corresponding to the same $d$.  We have
    $$
    p(t_{nkp},j)>\sup\tilde{V}-\frac{1}{p},\quad j\leq h_{np},\ j\in\k\bn.
    $$
    Since we assumed $h_{np}>n$, this implies
    $$
    p(t_{nkp},j)>\sup\tilde{V}-\frac{1}{p},\quad j\leq n,\ j\in\k\bn.
    $$
    Setting $\overline{t}_{n}=t_{nnn}$, we have 
    $\overline{t}_{n}>n$ hence $\overline{t}_{n}\to\infty$, and 
    $$
    \liminf_{n\in\bn}p(\overline{t}_{n},h)\geq\sup\tilde{V},
    \quad\forall h\in\k\bn.
    $$
    Then, by the proof of Proposition \ref{lim&ineq} $(ii)$, and Lemma
    \ref{Vtilde}, we have, for any $h\in\k\bn$,
    $$
    \limsup_{h\to+\infty} \limsup_{t\to+\infty}p(t,h)
    =\sup V\leq\sup\tilde{V}+\frac{2S}{\k}
    \leq \liminf_{n\in\bn}p(\overline{t}_{n},h)+\frac{2S}{\k}.
    $$
    Finally we observe that the function $\underline{g}(h)$ defined as 
    $\underline{g}(h)=\liminf_{n}g(\overline{t}_{n},h)$ is increasing, 
    therefore, if $\floor{\cdot}$ denotes the lower integer part, we get
    $\underline{g}(h)\geq\underline{g}(\floor{\frac{h}{\k}}\k)$, 
    from which
    $$
    \liminf_{h\to\infty}\frac{\underline{g}(h)}{h}
    \geq \liminf_{h\to\infty}\frac{\underline{g}
    \left(\floor{\frac{h}{\k}}\k\right)}{\floor{\frac{h}{\k}}\k}
    \frac{\floor{\frac{h}{\k}}\k}{h}=
    \liminf_{h\in\k\bn}\frac{\underline{g}(h)}{h},
    $$
    and the thesis follows.
\end{proof}

\end{document}